\documentclass{article}

\usepackage{amsmath,amsfonts,amssymb,a4,enumerate,epsfig,theorem,psfrag}
\usepackage[latin1]{inputenc}

\newcommand{\inv}{\textrm{\rm inv}}

\newcommand{\tub}{\textrm{\rm tub}}

\newcommand{\bbb}{{b}}

\newcommand{\dom}{{\operatorname{Dom}}}
\newcommand{\dist}{{\operatorname{dist}}}
\newcommand{\sign}{{\operatorname{sign}}}
\newcommand{\trace}{{\operatorname{tr}}}
\newcommand{\transpose}{{\ast}}

\newcommand{\NN}{{\mathbb{N}}}
\newcommand{\BB}{{\mathbb{B}}}

\newcommand{\RR}{{\mathbb{R}}}

\newcommand{\Mu}{{M}}
\newcommand{\Aa}{{\cal{A}}}
\newcommand{\Bb}{{\cal{B}}}
\newcommand{\Cc}{{\cal{C}}}
\newcommand{\Dd}{{\cal{D}}}
\newcommand{\Ee}{{\cal{E}}}

\newcommand{\Ii}{{\cal{I}}}
\newcommand{\Jj}{{\cal{J}}}
\newcommand{\Kk}{{\cal{K}}}

\newcommand{\Oo}{{\cal{O}}}

\newcommand{\calS}{{\cal{S}}}
\newcommand{\Tt}{{\cal{T}}}
\newcommand{\Uu}{{\cal{U}}}

\newcommand{\cX}{{\cal{X}}}

\newcommand{\cN}{{\cal{N}}}
\newcommand{\cM}{{\cal{M}}}
\newcommand{\Xx}{{\cal{X}}}

\newcommand{\ILa}{\Tt_\Lambda}
\newcommand{\OLa}{\Oo_\Lambda}
\newcommand{\OSa}{\Oo_\Sigma}

\newcommand{\ILai}{\Tt_{\Lambda_i}}

\newcommand{\JLa}{\bar\Jj_\Lambda}

\newcommand{\JLao}{\Jj_\Lambda}

\newcommand{\proof}{{\noindent\bf Proof: }}

\def\qed{\unskip\nobreak\hfil\penalty50\hskip1.75em\null\nobreak\hfil
$\blacksquare$ {\parfillskip=0pt \finalhyphendemerits=0 \par}\medbreak}

\newcommand\capsize{\relax}
\newcommand\nobf{\noindent\bf}

\newcommand\grad{\operatorname{grad}}
\newcommand\diag{\operatorname{diag}}

\newcommand\interior{\operatorname{int}}

\newtheorem{theo}{Theorem}
\newtheorem{lemma}{Lemma}[section]
\newtheorem{prop}[lemma]{Proposition}
\newtheorem{coro}[lemma]{Corollary}

\newtheorem{fact}[lemma]{Fact}

\title{Convergence rates to deflation of \\ simple shift strategies}
\author{Ricardo S. Leite, Nicolau C. Saldanha and Carlos Tomei}

\begin{document}

\maketitle

\begin{abstract}
The computation of eigenvalues of real symmetric tridiagonal
matrices frequently proceeds by a sequence of $QR$ steps with shifts.
We introduce \textit{simple shift strategies},
functions $\sigma$ satisfying natural conditions,
taking each $n \times n$ matrix $T$ to a real number $\sigma(T)$.
The strategy specifies the shift to be applied by the $QR$ step at $T$.
Rayleigh and Wilkinson's are examples of simple shift strategies.
We show that if $\sigma$ is continuous then there exist
initial conditions for which deflation does not occur,
i.e., subdiagonal entries do not tend to zero.
In case of deflation,
we consider the rate of convergence to zero of the $(n,n-1)$ entry:
for simple shift strategies this is always at least quadratic.
If the function $\sigma$ is smooth in a suitable region
and the spectrum of $T$ does not include three consecutive eigenvalues
in arithmetic progression then convergence is cubic.
This implies cubic convergence to deflation of Wilkinson's shift for generic spectra.
The study of the algorithm near deflation uses \textit{tubular coordinates},
under which $QR$ steps with shifts are given by a simple formula.
\end{abstract}

\medbreak


\smallbreak

{\noindent\bf MSC2010-class:} 65F15; 37N30.

\section{Introduction}

Let $\Tt$ be the vector space of real symmetric tridiagonal matrices. Among the
standard algorithms to compute eigenvalues of matrices in $\Tt$ are $QR$ steps
with different shift strategies: Rayleigh and Wilkinson are familiar examples
(excellent references are \cite{Wilkinson}, \cite{Demmel}, \cite{Parlett}).
In this paper, we consider a more general context: we define simple shift
strategies, which include the examples above and more, and discuss subtle aspects
of their asymptotic behavior.

More precisely, given a matrix $T \in \Tt$ and $s \in \RR$,
write $T - sI = Q R$, if possible, for an orthogonal matrix $Q$ and an upper triangular matrix $R$ with positive diagonal entries.
A \textit{shifted $QR$ step} is $\Phi(T,s) = Q^\transpose T Q$.
As is well known, shifted $QR$ steps preserve spectrum and shape.
For a real matrix with simple spectrum $\Lambda = \diag(\lambda_1,\ldots,\lambda_n)$,
let $\ILa \subset \Tt$ be the set of matrices similar to $\Lambda$.
A \textit{simple shift strategy} is a function $\sigma: \ILa \to \RR$
satisfying the following two properties.
\begin{enumerate}[(I)]
\item{For all $T \in \ILa$, $\sigma(E_n T E_n) =
\sigma(T)$, where $E_n = \diag(1, 1, \ldots, 1, -1)$.}
\item{There exists $C_\sigma > 0$ such that for all $T \in \ILa$
there is an eigenvalue $\lambda_i$ with $| \sigma(T) - \lambda_i | \le
C_\sigma |\bbb(T)|,$ where $\bbb(T) = (T)_{(n,n-1)}$.}
\end{enumerate}
This definition excludes algorithms which employ multi-shifts and extrapolation
techniques, but accomodates the usual Rayleigh and Wilkinson's strategies.

For technical reasons, we prefer the \textit{signed} variant
$\Phi_\star(T,s) = Q_\star^\transpose T Q_\star$,
where now $T - sI = Q_\star R_\star$,
the orthogonal matrix $Q_\star$ has positive determinant and
only the first $n-1$ diagonal entries of the upper triangular matrix $R_\star$
are required to be positive. It is easy to see that  either
$\Phi(T,s)=\Phi_\star(T,s)$ or $\Phi(T,s)= E_n \Phi_\star(T,s) E_n$
(notice that $E_n T E_n$ is obtained from $T \in \Tt$
by changing the signs of the entries $(n,n-1)$ and $(n-1,n)$).
As we shall see, the signed step is smoothly defined on a larger domain,
and convergence issues for both kinds of step iterations are essentially equivalent.


Simple shift strategies prescribe shifts: set
$F_s(T) = \Phi_\star(T,s)$ and $F_\sigma(T)= F_{\sigma(T)}(T)$.
Algorithms iterate $F_\sigma$ aiming at deflation, i.e.,
obtaining a matrix $T \in \ILa$ with small $|\bbb(T)|$.
The \textit{deflation set} (resp. \textit{neighborhood}) is the set
$\Dd_{\Lambda,0} \subset \ILa$
(resp. $\Dd_{\Lambda,\epsilon} \subset \ILa$)
consisting of matrices $T$ for which
$\bbb(T)=0$ (resp. $|\bbb(T)| \le \epsilon$).
A simple shift strategy
$\sigma$ is \textit{deflationary} if for any $T \in \ILa$ and any $\epsilon >0$
there exists $k$ for which $F_\sigma^k(T) \in \Dd_{\Lambda,\epsilon}$.
As is well known, Rayleigh's strategy is not deflationary and Wilkinson's is:
our first result provides a context for these facts.

\begin{theo}\label{theo:connect} A continuous shift strategy $\sigma: \ILa \to
\RR$ is not deflationary.
\end{theo}

We are thus led to consider the \textit{singular support}
$\calS_\sigma \subset \ILa$ of a shift strategy $\sigma$,
i.e., the minimal closed subset of $\ILa$
on whose complement $\sigma$ is smooth.
For Rayleigh $\calS_\sigma$ is empty;
for Wilkinson's strategy it consists of the
matrices $T \in \ILa$ with $(T)_{n-1,n-1} = (T)_{n,n}$.

Numerical evidence brings up the question of whether
the rate of convergence to zero of the sequence $\bbb(F_\sigma^k(T))$ is cubic,
in the sense that there is a constant $C$ such that
$|\bbb(F_\sigma^{k+1}(T))| \le C |\bbb(F_\sigma^k(T))|^3$ for large $k$.
It turns out that, for any shift strategy $\sigma$,
in an appropriate neighborhood of the deflation set $\Dd_{\Lambda,0}$,
each iteration of $F_\sigma$ squeezes the $(n,n-1)$ entry quadratically.
Away from the singular support $\calS_\sigma$, squeezing is cubic.

\begin{theo}
\label{theo:squeeze}
For $\epsilon > 0$ small enough,
the deflation neighborhood $\Dd_{\Lambda,\epsilon}$
is invariant under $F_\sigma$.
There exists $C > 0$ such that,
for all $T \in \Dd_{\Lambda,\epsilon}$,
$|\bbb(F_\sigma(T))| \le C |\bbb(T)|^2$.
Also, given a compact set $\Kk\subset \Dd_{\Lambda,\epsilon}$
disjoint from $\calS_\sigma \cap \Dd_{\Lambda,0}$,
there exists $C_\Kk > 0$ such that,
for all $T \in \Kk$,
$|\bbb(F_\sigma(T))| \le C_{\Kk} |\bbb(T)|^3$.
\end{theo}

Cubic convergence does not hold in general for Wilkinson's strategy.
In \cite{LST2}, for $\Lambda= \diag(-1,0,1)$,
we construct a Cantor-like set $\cX \subset \ILa$
of unreduced initial conditions
for which the rate of convergence is strictly quadratic.
Sequences starting at $\cX$ converge to a reduced
matrix which is not diagonal.
For Rayleigh's shift, on the other hand,
convergence is always cubic within invariant deflation neighborhoods.

A matrix $T \in \Tt$ with simple spectrum is {\it a.p.\ free} if it does not
have three eigenvalues in arithmetic progression and {\it a.p.}\ otherwise.
For a.p.\ free spectra, the situation is very nice:
cubic convergence is essentially uniform on $\ILa$.


\begin{theo}
\label{theo:big}
Let $\Lambda$ be an a.p.\ free matrix and
$\sigma$ a shift strategy
for which diagonal matrices do not belong to $\calS_\sigma$.
Then there exist $\epsilon > 0$, $C >0$ and $K > 0$ such that
the deflation neighborhood $\Dd_{\Lambda,\epsilon}$ is invariant under $F_\sigma$.
Also, for any $T \in \Dd_{\Lambda,\epsilon}$,
the sequence $(F_\sigma^k(T))$ converges to a diagonal matrix
and the set of positive integers $k$ for which
$|\bbb(F_\sigma^{k+1}(T))| > C |\bbb(F_\sigma^{k}(T))|^3$
has at most $K$ elements.
\end{theo}

Still, the finite set of points in which the cubic estimate does not
hold may occur arbitrarily late along the sequence $(F_\sigma^{k}(T))$.

An a.p.\ matrix is {\it strong a.p.}\ if it contains three consecutive
eigenvalues in arithmetic progression and {\it weak a.p.}\ otherwise.
Under very mild additional hypothesis, $\bbb(T)$ converges to zero at a cubic
rate also for weak a.p.\ matrices.
Let $\Cc_{\Lambda,0} \subset \ILa$ be the set of matrices $T$
for which $(T)_{n,n-1} = (T)_{n-1,n-2} = 0$.

\begin{theo}
\label{theo:bigg}
Let $\Lambda$ be a weak a.p.\ matrix and
$\sigma: \ILa \to \RR$ a shift strategy
for which $\Cc_{\Lambda,0}$ and $\calS_\sigma$ are disjoint.
Then there exists $\epsilon > 0$ such that
the deflation neighborhood  $\Dd_{\Lambda,\epsilon}$ is invariant under $F_\sigma$
and, for all unreduced $T \in \Dd_{\Lambda,\epsilon}$,
the sequence $(\bbb(F_\sigma^k(T)))$ converges to zero
at a rate which is at least cubic.
More precisely,
for each unreduced $T \in \Dd_{\Lambda,\epsilon}$
there exist $C_T, K_T > 0$ such that,
for all $k > K_T$,
we have $|\bbb(F_\sigma^{k+1}(T))| \le C_T |\bbb(F_\sigma^k(T))|^3$.
\end{theo}

In particular, the convergence of Wilkinson's strategy is cubic for weak a.p.\ matrices.
However, uniformity in the sense of Theorem \ref{theo:big} is not guaranteed and
the constants $C_T$ and $K_T$ depend on $T$.
As in the case of the spectrum $\{-1,0,1\}$,
we conjecture that if $\Lambda$ is strong a.p.\ then there exists $\cX \subset \ILa$
of Hausdorff codimension 1 of initial conditions $T$ for which the rate of
convergence is strictly quadratic.

\smallskip

The proofs of the above results depend on few basic ideas.
Signed shifted steps $F_s(T)$ are shown to be well defined for unreduced matrices
and in an open neighborhood of the deflation set $\Dd_{\Lambda,0}$.
The compact set (manifold) $\Dd_{\Lambda,0}$ splits into connected
components $\Dd^i_{\Lambda,0}$
consisting of matrices $T \in \ILa$ with $(T)_{n,n}=\lambda_i$.
For small $\epsilon$, the deflation neighborhood splits into components $\Dd_{\Lambda,\epsilon}^i$ which are thickenings of $\Dd^i_{\Lambda,0}$.

\textit{Tubular coordinates} provide a good understanding both of
tubular neighborhoods $\Dd_{\Lambda,\epsilon}^i$
and of shifted $QR$ steps within these sets.
A previous unpublished version of this paper (\cite{LST3})
uses instead \textit{bidiagonal coordinates}, defined in \cite{LST1},
to prove some of the results presented here for Wilkinson's shift;
these coordinates are also used in \cite{LST1}
to prove the cubic convergence of Rayleigh's shift.
Bidiagonal coordinates consist of very explicit charts on the manifold $\ILa$.
On both coordinate systems, shifted $QR$ steps are very simple.

Sufficiently thin deflation neighborhoods are invariant under steps $F_\sigma$.
Theorem \ref{theo:connect} then becomes a connectivity argument:
in a nutshell, there must be separatrices
in order to permit deflation to different components $\Dd^i_{\Lambda,0}$.

Steps $F_\sigma$ are smooth whenever the shift strategy is, i.e.,
for $T \in \Dd_{\Lambda,\epsilon}^i \smallsetminus \calS_\sigma$.
At matrices $T_0 \in \Dd^i_{\Lambda,0}$ on which $F_\sigma$ is smooth,
the map $T \mapsto \bbb(F_\sigma(T))$ has zero gradient.
The symmetry of the shift strategy (condition (I))
yields a cubic Taylor expansion and therefore an estimate
$|\bbb(F_\sigma(T))| \le C |\bbb(T)|^3$, settling Theorem \ref{theo:squeeze}.

\textit{Height functions} $H: \Dd_{\Lambda,\epsilon}^i \to \RR$
are used for further study of the sequence $(F_\sigma^k(T))$.
More precisely, for steps $s$ near $\lambda_i$,
$H_i(F_s(T))> H_i(T)$ provided $T \in \Dd^i_{\Lambda,\epsilon}$ is not diagonal.
A compactness argument then bounds the number of iterations
for which $F_\sigma^k(T)$ stays close to the singular support ${\calS}_\sigma$:
this is essentially Theorem \ref{theo:big}.

For a.p. spectra, the situation is subtler.
As numerical analysts know,
shift strategies usually define sequences of matrices which,
asymptotically, not only isolate an eigenvalue at the $(n,n)$ position but also isolate,
at a slower rate, a second eigenvalue at the $(n-1,n-1)$ position.
This does not happen for the example in \cite{LST2} where
$(F_\sigma^k(T))_{n,n}$ tends to the center
of a three-term arithmetic progression of eigenvalues and
$(F_\sigma^k(T))_{n-1,n-2}$ stays bounded away from $0$.
On the other hand, Theorem \ref{theo:bigg} tells us that
the weak a.p.\ hypothesis
together with an appropriate smoothness condition
guarantee cubic convergence.

\smallskip

In Section 2 we list the basic properties of the signed shifted $QR$ step.
Simple shift strategies are introduced in Section 3,
and the standard examples are shown to satisfy the definition.
We define the deflation set $\Dd_{\Lambda,0}$ and neighborhood $\Dd_{\Lambda,\epsilon}$
in Section 4 and then set up tubular coordinates.
The local theory of steps $F_s$ near $\Dd_{\Lambda,0}$
and the proof Theorem \ref{theo:squeeze} are presented in Section 5.
Section 6 is dedicated to Theorem \ref{theo:connect}.
In Section 7 we construct the height functions $H$ and then prove Theorem \ref{theo:big}.
The convergence properties of a.p.\ matrices
in Theorem \ref{theo:bigg} are proved in Section 8.
We finally present in Section 9 two counterexamples
to natural but incorrect strengthenings of Theorems \ref{theo:big} and \ref{theo:bigg}.


The authors are very grateful for the abundant contributions 
of several readers of this work and its previous versions.
The authors acknowledge support from CNPq, CAPES, IM-AGIMB and FAPERJ.


\section{$QR$ iteration with shift and a variation}

For a matrix $M$, the $QR$ factorization is $M = QR$
for an orthogonal matrix $Q$ and
an upper triangular matrix $R$ with positive diagonal.
As usual, let $SO(n)$ denote the set of orthogonal matrices
with determinant equal to $1$.
The {\it $Q_\star R_\star$ factorization}, instead,
is $M = Q_\star R_\star$, for $Q_\star \in SO(n)$ and
$R_\star$ an upper triangular matrix
with $(R_\star)_{i,i} > 0$, $i=1,\ldots,n-1$.
A real $n \times n$ matrix $M$ is {\it almost invertible}
if its first $n-1$ columns are linearly independent:
notice that almost invertible matrices are dense
within $n \times n$ matrices and form an open set.
The diagonal matrix $E_n$ is such that $(E_n)_{i,i}$
is $1$ for $i < n$ and $-1$ for $i=n$.

\begin{prop}
\label{prop:existence}
An almost invertible real matrix $M$ admits
a unique $Q_\star R_\star$ factorization,
with $Q_\star$ and $R_\star$ depending smoothly
on $M$.
If $M$ is invertible, it admits unique (smooth)
factorizations $M = QR = Q_\star R_\star$.
If $\det M >0$, the factorizations are equal, i.e.,
$Q= Q_\star$ and  $R = R_\star$.
If $\det M <0$, $Q= Q_\star E_n$ and $R = E_n R_\star$.
If $\det M =0$, $(R_\star)_{n,n}=0$.
\end{prop}


\proof
Let $M$ be almost invertible. Applying Gram-Schmidt with positive
normalizations on its first $n-1$ columns we obtain the first $n-1$ columns of
both $Q$ and $R$, as well as those of $Q_\star$ and $R_\star$. The last column
$v = Q_\star e_n$ of $Q_\star$ is already well defined, by orthonormality and
the fact that $\det Q_\star = 1$. Now, set $R_\star = M (Q_\star)^\ast$. The
positivity of $R_{n,n}$ specifies whether the last column of $Q$ is $v$ or
$-v$. Smoothness
is clear by construction.

If $M$ is invertible, $\det M = \det Q_\star  \det R_\star$ implies that the
last diagonal entry of $R_\star$ has the same sign of $\det M$: the relations
between the factorizations then follow. If $M$ is not invertible, the relation
among determinants implies $(R_\star)_{n,n}=0$.
\qed

Let $\Tt$ denote the real vector space of $n \times n$ real, symmetric,
tridiagonal matrices endowed with the norm $\|T\|^2 = \trace(T^2)$.
For $T \in \Tt$, the \textit{subdiagonal entries}
of $T$ are $(T)_{i+1,i}$ for $i=1,\ldots,n-1$.
The lowest subdiagonal entry of $T$ is $\bbb(T)= (T)_{n,n-1}$.
If all subdiagonal entries of $T$ are nonzero, $T$ is
an \textit{unreduced} matrix; otherwise, $T$ is {\it reduced}.
Notice that an unreduced tridiagonal matrix is almost invertible:
indeed, the block formed by
rows $2,\ldots,n$ and columns $1,\ldots,n-1$ is a an upper triangular matrix
with nonzero diagonal entries, and therefore, invertible.

We consider the {\it shifted $QR$ step} and its {\it signed} counterpart,
\[ \Phi(T, s) = Q^\ast T Q, \quad \Phi_\star(T, s) = Q_\star^\ast T Q_\star, \]
where $T - s I = QR$ and $T - s I = Q_\star R_\star$.
The pair $(T,s) \in \Tt \times \RR$ belongs to the natural (open, dense)
domains $\dom(\Phi)$ and $\dom(\Phi_\star)$
if $T- sI$ is invertible and almost invertible, respectively:
clearly, the functions $\Phi$ and $\Phi_\star$ are smooth
in their domains.

\begin{lemma}
\label{lemma:basic}
For $(T,s) \in \dom(\Phi)$ (resp.  $\dom(\Phi_\star)$),
we have $\Phi(T,s) \in \Tt$ (resp. $\Phi_\star(T,s) \in \Tt$).
The spectra of $T$, $\Phi(T,s)$ and $\Phi_\star(T,s)$ are equal.
In the appropriate domains,
for $T - sI = QR = Q_\star R_\star$ and $i = 1, 2, \ldots, n-1$,
\[ (\Phi(T,s))_{i+1,i} =
\frac{(R)_{i+1,i+1}}{(R)_{i,i}} \; (T)_{i+1,i}, \quad
(\Phi_\star(T,s))_{i+1,i} =
\frac{(R_\star)_{i+1,i+1}}{(R_\star)_{i,i}} \; (T)_{i+1,i}. \]
Thus, the top $n-2$ subdiagonal entries of $T$,
$\Phi(T,s)$ and $\Phi_\star(T,s)$ have the same sign;
also, $\sign\, (T)_{n,n-1}  = \sign\, (\Phi(T,s))_{n,n-1}$.

\end{lemma}

\smallskip

\proof
We prove the statements for $\Phi_\star$; the others are then easy.

For a pair $(T,s) \in \dom (\Phi) \subset \dom (\Phi_\star)$,
there are two expressions for $\Phi_\star(T,s)$:
\[ \Phi_\star(T,s) = Q_\star^\ast T Q_\star = R_\star T R_\star^{-1},
\quad \hbox{where}\quad T - sI = Q_\star R_\star. \]
From the first equality, $\Phi_\star(T,s)$ is symmetric and
from the second, $\Phi_\star(T,s)$ is an upper Hessenberg matrix
so that $\Phi_\star(T,s) \in \Tt$ is similar to $T$.
More generally, for $(T,s) \in \dom (\Phi_\star)$
we still have
\[ \Phi_\star(T,s) = Q^\ast_\star T Q_\star , \quad
\Phi_\star(T,s) R_\star = R_\star T \]
and therefore $\Phi_\star(T,s) \in \Tt$ is similar to $T$.
Compute the $(i+1,i)$ entry of the second equation above to obtain
$(\Phi_\star(T,s))_{i+1,i} \; (R_\star)_{i,i} =
(R_\star)_{i+1,i+1} \; (T)_{i+1,i}$,
completing the proof.
\qed

The following result describes the behavior of $\Phi_\star$
at points not in $\dom(\Phi)$,
which will play an important role throughout the paper.

\begin{lemma}
\label{lemma:sislambda}
If $(T,s) \in \dom(\Phi_\star) \smallsetminus \dom(\Phi)$ then
\[ \bbb(\Phi_\star(T,s)) = (\Phi_\star(T,s))_{n,n-1} = 0, \quad
(\Phi_\star(T,s))_{n,n} = s. \]
At a point $(T,s) \in \dom(\Phi_\star)$ with
$\bbb(T) = 0$ and $s = (T)_{n,n}$ we
have $\grad(\bbb \circ \Phi_\star) = 0$.
\end{lemma}

\proof
Since $T - sI = Q_\star R_\star = R_\star^\ast Q_\star^\ast$ is not invertible
then $(R_\star)_{n,n} = 0$ and therefore $R_\star^\ast e_n = 0$.
Thus $v = (Q_\star^\ast)^{-1} e_n = Q e_n$ satisfies $(T - sI)v = 0$.
We then have
$\Phi_\star(T,s) e_n = Q^\ast T Q e_n = Q^\ast T v = Q^\ast (sv) = se_n$,
proving the first claim.
For the second claim, since $T - sI$ is almost invertible,
$(R_\star)_{i,i} > 0$ for $i < n$.
From the previous lemma,
\[ (\bbb \circ \Phi_\star)(T,s) =
\frac{(R_\star)_{n,n}}{(R_\star)_{n-1,n-1}} \; \bbb(T)\,; \]
if $\bbb(T) = 0$ and $s = (T)_{n,n}$ then
$(R_\star)_{n,n} = 0$ and
$\bbb \circ \Phi_\star$ is a product of two smooth functions,
both zero, yielding $\grad(\bbb \circ \Phi_\star) = 0$.
\qed

Recall that symmetrically changing the signs of subdiagonal entries
does not change the spectrum of a matrix in $\Tt$.
Let $\Ee$ denote the set of \textit{signed diagonal matrices}
with $\pm 1$ along the diagonal entries;
in particular, $E_n \in \Ee$.
The operation of changing subdiagonal signs, i.e.,
of conjugation by some $E \in \Ee$,
behaves well with respect to $\Phi$ and $\Phi_\star$.

\begin{lemma}
\label{lemma:ETE}
Let $E \in \Ee$.
The domains $\dom(\Phi)$ and $\dom(\Phi_\star)$
are invariant under conjugation by $E$ and
\[ \Phi(ETE,s)= E \Phi(T,s) E, \quad
\Phi_\star(ETE,s)= E \Phi_\star(T,s ) E. \]
If $\det (T - s I) > 0 $ then $\Phi(T,s)=\Phi_\star(T,s)$;
if  $\det (T- s I) < 0 $, $\Phi(T,s)= E_n \Phi_\star(T,s) E_n$;
if  $\det (T- s I) = 0 $ and $(T,s) \in \dom(\Phi_\star)$,
then $\bbb(\Phi_\star(T,s))=0$.
\end{lemma}


\proof
For $(T,s) \in \dom(\Phi)$,
the matrices $T - sI$ and $E(T - sI)E$ are both invertible.
The $QR$ factorization $T -sI = Q R$ yields
$ETE - E(sI)E = (EQE)(ERE)$,
preserving the positivity of the diagonal entries of the triangular part, so
\[ \Phi(ETE,s) = (EQE)^\ast ETE (EQE) = E Q^\ast T Q E = E \Phi(T,s) E.\]
The argument is similar for $\Phi_\star$.
The claims for $T-sI$ invertible follow from the relation
between $Q$ and $Q_\star$ in Proposition \ref{prop:existence};
the case $\det(T - sI) = 0$ is a repetition of Lemma \ref{lemma:sislambda}.
\qed

We are only interested in the case when the spectrum of $T$ is simple,
since a double eigenvalue implies reducibility.
Let $\Lambda$ be a real diagonal matrix
with simple eigenvalues $\lambda_1 < \cdots < \lambda_n$.
Define the \textit{isospectral manifold}
\[ \ILa = \{ Q^\ast \Lambda Q, Q \in SO(n)\} \cap \Tt, \]
the set of matrices in $\Tt$ similar to $\Lambda$.
The set $\ILa \subset \Tt$ is a real smooth
manifold (\cite{Tomei}; \cite{LST1} describes an explicit atlas of $\ILa$).
Since either version of shifted $QR$ step preserves spectrum,
restriction defines smooth
maps
$\Phi: (\ILa \times \RR) \cap \dom(\Phi) \to \ILa$ and
$\Phi_\star: (\ILa \times \RR) \cap \dom(\Phi_\star) \to \ILa$.

Still in $\ILa$, it is convenient to consider
the \textit{step} $F_s(T) = \Phi_\star(T,s)$.
For $s$ not an eigenvalue of $\Lambda$, the domain of $F_s$ is $\ILa$.
The natural domain for $F_{\lambda_i}$ instead
is the \textit{deflation domain} $\Dd_\Lambda^i$,
the open dense subset of $\ILa$ of matrices $T$
for which $T - \lambda_i I$ is almost invertible.
In other words, $T \in \Dd_\Lambda^i$ if and only if
$\lambda_i$ is an eigenvalue of the lowest irreducible block of $T$.

The definition of the step $F_s$ differs from the usual one
in that we use $\Phi_\star$ instead of $\Phi$.
Given Lemma \ref{lemma:ETE},
considerations about deflation are unaffected
and our choice has the advantage of being smooth
(and well defined) in $\Dd_{\Lambda}^i$.

The ($i$-th) \textit{deflation set} is
\[ \Dd_{\Lambda,0}^i =
\left\{ T \in \ILa\;|\; \bbb(T)=0, (T)_{n,n} = \lambda_i \right\}. \]
Since the spectrum of $\Lambda$ is simple,
$\Dd_{\Lambda,0}^i \subset \Dd_\Lambda^i$.
Also, if $i \ne j$ then $\Dd_\Lambda^i \cap \Dd_{\Lambda,0}^j = \emptyset$.
We saw in Lemma \ref{lemma:sislambda} that
when the shift is taken to be an eigenvalue, a single step deflates a matrix,
i.e., that the image of $F_{\lambda_i}$ is contained in $\Dd^i_{\Lambda,0}$:
we shall see in Proposition \ref{prop:Fs} that this image
is in fact equal to $\Dd^i_{\Lambda,0}$.

\section{Simple shift strategies}

Quoting Parlett \cite{Parlett}, there are shifts for all seasons.
The point of using a shift strategy is to accelerate deflation,
ideally by choosing $s$ near an eigenvalue of $T$.
A \textit{simple shift strategy} is a function
$\sigma: \ILa \to \RR$ such that:
\begin{enumerate}[(I)]
\item{for all $T \in \ILa$, $\sigma(E_n T E_n) = \sigma(T)$;}
\item{there exists $C_\sigma > 0$ such that
for all $T \in \ILa$ there is an eigenvalue $\lambda_i$ with
$| \sigma(T) - \lambda_i | \le C_\sigma |\bbb(T)|.$}
\end{enumerate}
In particular, if $T \in \Dd_{\Lambda,0}^i$ then $\sigma(T) = \lambda_i$.
The \textit{step} associated with a (simple) shift strategy $\sigma$
is $F_\sigma$, defined by $F_\sigma(T) = F_{\sigma(T)}(T)$.
The natural domain for $F_\sigma$ is the set of matrices $T$
for which $T - \sigma(T) I$ is almost invertible.
From Section 2, it includes all unreduced matrices
and open neighborhoods of each deflation set $\Dd^i_{\Lambda,0}$.
We shall also see in Section 6 that it contains
a dense open subset $\Uu_{\Lambda,\epsilon}$ of $\ILa$ invariant under $F_\sigma$.  
A more careful description of this domain will not be needed.


We leave to the reader the verification that
\textit{Rayleigh's shift} $\rho(T) = (T)_{n,n}$ is a simple shift strategy.
Denote the bottom $2 \times 2$ diagonal principal minor
of a matrix $T \in \Tt$ by $\hat{T}$:
\textit{Wilkinson's shift} $\omega(T)$ is the eigenvalue
of $\hat T$ closer to $(T)_{n,n}$
(in case of draw, take the smallest eigenvalue).

\begin{lemma}
\label{lemma:2sqrt2}
The function $\omega$ is a simple shift
strategy with $C_\omega = 2\sqrt{2}$.
\end{lemma}

\proof
Condition (I) follows from the fact that changing signs of off-diagonal
entries of a $2 \times 2$ matrix does not change its spectrum.
For (II), apply the Wielandt-Hoffman theorem
to the $2 \times 2$ trailing principal minors of $T$ and
$S = T - \bbb(T) B$ to deduce that
$|(T)_{n,n} - \omega(T)| \le \sqrt{2} \; \bbb(T)$.
Again from Wielandt-Hoffman,
now on $S$ as in Proposition \ref{prop:deflation},
$|(T)_{n,n} - \lambda_i| \le \sqrt{2} \; \bbb(T)$ and
$| \omega(T) - \lambda_i| \le 2\sqrt{2} \; \bbb(T)$.
\qed

Another example of shift strategy,
the \textit{mixed Wilkinson-Rayleigh strategy},
uses Wilkinson's shift unless the matrix is already near deflation,
in which case we use Rayleigh's:
\[ \sigma(T) = \begin{cases}
\rho(T),&|(T)_{n,n-1}| < \epsilon, \\
\omega(T),&|(T)_{n,n-1}| \ge \epsilon;
\end{cases} \]
here $\epsilon > 0$ is a small constant.

Shift strategies are not required to be continuous and
$\omega$ is definitely not.
For a shift strategy $\sigma$,
let $\calS_\sigma \subset \ILa$ be the \textit{singular support} of $\sigma$,
i.e., a minimal closed set on whose complement $\sigma$ is smooth.
For example, $\calS_\omega$ is the set of matrices $T \in \ILa$
for which the two eigenvalues $\omega_-(T)$ and $\omega_+(T)$ of $\hat T$
are equidistant from $(T)_{n,n}$, or,
equivalently, for which $(T)_{n,n} = (T)_{n-1,n-1}$.
The set $\calS_\sigma$ will play an important role later.

We consider the phase portrait of $F_\omega$ for $3 \times 3$ matrices.
In this case,
the reader may check that the domain of $F_\omega$ is the full set $\ILa$.
Let $\JLao \subset \ILa$ be set of \textit{Jacobi matrices}
similar to $\Lambda$, i.e., matrices $T \in \ILa$
with strictly positive subdiagonal entries.
It is known (\cite{DNT}, \cite{Moerbeke}) that
the closure $\JLa \subset \ILa$ is diffeomorphic to a hexagon.
The set $\JLa$ is not invariant under $F_\omega$
but we may define $\tilde F_\omega(T)$ with $\tilde F_\omega: \JLa \to \JLa$
by dropping signs of subdiagonal entries of $F_\omega(T)$.
As discussed above, this rather standard procedure is mostly harmless.

\begin{figure}[ht]
\begin{center}
\psfrag{Yy}{$\calS_\omega$}
\epsfig{height=28mm,file=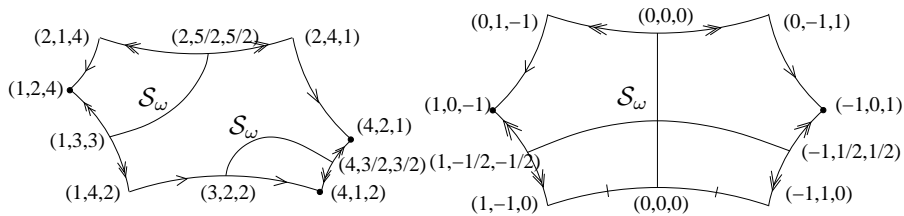}
\end{center}
\caption{\capsize The phase space of Wilkinson's step for $n=3$.}
\label{fig:wil}
\end{figure}

Two examples of $\tilde F_\omega$ are given in Figure \ref{fig:wil},
which represent $\JLa$ for the $\Lambda= \diag(1,2,4)$ on the left
and $\Lambda=\diag(-1,0,1)$ on the right.
The vertices are the six diagonal matrices
similar to $\Lambda$ and the edges consist of reduced matrices.
Labels indicate the diagonal entries of the corresponding matrices.
Three edges form $\Dd^i_{\Lambda,0} \cap \JLa$:
they alternate, starting from the bottom horizontal edge on both hexagons.
The set $\calS \cap \JLa$ is indicated in both cases.


Vertices are fixed points of $\tilde F_\omega$ and boundary edges are invariant sets.
A simple arrow indicates the motion of the points $\tilde F_\omega^k(T)$ along the edge.
Points $T$ on an arc with a double arrow are taken to a diagonal matrix in a single step:
the arc points to $\tilde F_\omega(T)$.
Arcs marked with a transversal segment consist of fixed points of $\tilde F_\omega$.

Points on both sides of $\calS_\omega$ are taken far apart:
there is essentially a jump discontinuity along $\calS_\omega$.
From Theorem \ref{theo:squeeze},
the decay of the bottom subdiagonal entry under Wilkinson's step
away from $\calS_\omega \cap \Dd_{\Lambda,0}$ is cubic.
As discussed in \cite{LST2}, near $\calS_\omega \cap \Dd_{\Lambda,0}$
this decay is quadratic, but not cubic.
For the left hexagon, cubic convergence occurs in the long run
because the sequence $\tilde F_\omega^k(T)$ stays close to this intersection
only for a few values of $k$,
illustrating Theorem \ref{theo:big}.

In the case $\Lambda= \diag(-1,0,1)$, the bottom edge consists of fixed points.
This gives rise to a special asymptotic behavior (\cite{LST2}):
the (fixed) point labeled by $(0,0,0)$ is actually
the limit of a collection of sequences $(\tilde F_\omega^k(T))$
for which the convergence is strictly quadratic.


\section{Tubular coordinates}

We collect a few basic facts about shifted $QR$ steps.

\begin{prop}
\label{prop:Fs}
If $s$ is not an eigenvalue of $\Lambda$,
the map $F_s: \ILa \to \ILa$ is a diffeomorphism.
The image of $F_{\lambda_i}: \Dd_\Lambda^i \to \ILa$ is $\Dd_{\Lambda,0}^i$.
The restriction
$F_{\lambda_i}|_{\Dd^i_{\Lambda,0}}: \Dd^i_{\Lambda,0} \to \Dd^i_{\Lambda,0}$
is a diffeomorphism.
\end{prop}

\proof
If $s$ is not an eigenvalue,
compute $F_s^{-1}(T)$ by factoring $T - sI$ as $RQ$,
$R$ upper triangular with the first $n-1$ diagonal entries positive
and $Q \in SO(n)$:
we claim that $F_s(T_0) = T$ for $T_0 = QR + sI$, proving
that $F_s$ is a diffeomorphism.
Indeed, $QR = T_0 - sI$ is a $Q_\star R_\star$
factorization and thus $F_s(T_0) = Q^\ast T_0 Q = T$.

From the last sentence of Section 2,
the image of $F_{\lambda_i}$ is contained in
$\Dd^i_{\Lambda,0} \subset \Dd^i_{\Lambda}$.
The fact that the restriction of $F_{\lambda_i}$ to $\Dd^i_{\Lambda,0}$ is
a diffeomorphism is proved as in the previous paragraph.
\qed

Commutativity of steps is well known
and related to the complete integrability
of the interpolating Toda flows
(\cite{Flaschka}, \cite{LT}, \cite{Moser}, \cite{Parlett}).
For the reader's convenience we provide a proof.

\begin{prop}
\label{prop:commute}
Steps commute:
$F_{s_0} \circ F_{s_1} = F_{s_1} \circ F_{s_0}$
in the appropriate domains.
\end{prop}

The domain of $F_{s_0} \circ F_{s_1} = F_{s_1} \circ F_{s_0}$ is $\ILa$ if
neither $s_0$ nor $s_1$ is an eigenvalue,
$\Dd^i_\Lambda$ if $s_0 = \lambda_i$
and $s_1$ is not an eigenvalue (or vice-versa) and
the empty set in the rather pointless case
$s_0 = \lambda_i$, $s_1 = \lambda_j$, $i \ne j$.

\smallskip

\proof
We prove commutativity only when $s_0$ and $s_1$ are not eigenvalues;
the other cases follow easily.
Consider $Q_\star R_\star$ factorizations
\begin{gather*}
T - s_0 I = Q_0 R_0, \quad T - s_1 I = Q_1 R_1, \\
(T - s_0 I)(T - s_1 I) = (T - s_1 I)(T - s_0 I) = Q_2 R_2.
\end{gather*}
For $F_{s_0}(T) - s_1 = Q_0^\ast(T - s_1)Q_0 = Q_3 R_3$,
we have $F_{s_1}(F_{s_0}(T)) = Q_3^\ast F_{s_0}(T) Q_3 =
Q_3^\ast Q_0^\ast T Q_0 Q_3$.
Thus
\[ Q_0^\ast(T - s_1)Q_0R_0 = Q_0^\ast (T - s_1 I)(T - s_0 I) =
Q_0^\ast Q_2 R_2 = Q_3 R_3 R_0 \]
and therefore $Q_0^\ast Q_2 = Q_3$ and
$F_{s_0}(F_{s_1}(T)) = Q_2^\ast T Q_2$.
\qed

Recall that a map $\Pi: X \to Y \subset X$ is a \textit{projection}
if $\Pi(X) = Y$ and $\Pi \circ \Pi = \Pi$.
The map $F_{\lambda_i}: \Dd_\Lambda^i \to \Dd^i_{\Lambda,0}$ 
is not a projection but can be used to define one:
the \textit{canonical projection} 
$\Pi_i: \Dd_\Lambda^i \to \Dd^i_{\Lambda,0}$,
\[ \Pi_i(T) =
(F_{\lambda_i}|_{\Dd^i_{\Lambda,0}})^{-1}(F_{\lambda_i}(T)). \]

\begin{prop}
\label{prop:pia}
The map $\Pi_i$ is indeed a
smooth projection which commutes with steps.
More precisely,
$\Pi_i(F_s(T)) = F_s(\Pi_i(T))$ provided
$s$ is not an eigenvalue of $\Lambda$ different from $\lambda_i$.
\end{prop}

\proof
The map $\Pi_i$ is clearly smooth and, for $T \in
\Dd^i_{\Lambda,0}$, we have \[ \Pi_i(T) =
(F_{\lambda_i}|_{\Dd^i_{\Lambda,0}})^{-1}(F_{\lambda_i}(T)) = T, \] proving
that $\Pi_i$ is a projection.  Commutativity follows from
Proposition \ref{prop:commute}.
\qed

For a diagonal matrix $\Lambda$ with simple spectrum and $\epsilon > 0$,
the \textit{deflation neighborhood}
$\Dd_{\Lambda,\epsilon}\subset \ILa$
is the closed set of matrices $T \in \ILa$ with $|\bbb(T)| \le \epsilon$.
This notation is consistent with $\Dd_{\Lambda,0}$
for the deflation set.
As we shall see in Propositions \ref{prop:deflation} and \ref{prop:epsilonF},
for sufficiently small $\epsilon > 0$ the set $\Dd_{\Lambda,\epsilon}$
has connected components $\Dd^i_{\Lambda,\epsilon} \subset \Dd^i_\Lambda$,
$\Dd^i_{\Lambda,\epsilon} \supset \Dd^i_{\Lambda,0}$,
which are invariant under steps $F_s$ for shifts $s$ near $\lambda_i$, i.e.,
$F_s(\Dd^i_{\Lambda,\epsilon}) \subset \Dd^i_{\Lambda,\epsilon}$.
The sets $\Dd^i_{\Lambda,\epsilon}$ are therefore
also invariant under $F_\sigma$.

Denote the distance between a matrix $T$
and a compact set of matrices $\cN$ by
$\dist(T,\cN) = \min_{S \in \cN} \| T - S \|$.
Let $\gamma = \min_{i \ne j} |\lambda_i - \lambda_j|$ be
the {\it spectral gap} of $\Lambda$
and $B = e_n e_{n-1}^\ast + e_{n-1} e_n^\ast$.

Recall that if $\cN$ is a submanifold of codimension $k$ of $\cM$
then a \textit{closed tubular neighborhood} of $\cN$
consists of a closed neighborhood $\cN_\epsilon$ of $\cN$
and a diffeomorphism
$\zeta: \cN_\epsilon \to \cN \times \BB^k_\epsilon$ with
$\zeta(x) = (x,0)$ for $x \in \cN$
(here $\BB^k_\epsilon \subset \RR^k$ is the closed ball
of radius $\epsilon$ around the origin).
Given $x \in \cN$, the preimage
$\zeta^{-1}(\{x\} \times \BB^k_\epsilon)$
is a manifold with boundary of dimension $k$,
the \textit{fiber} through $x$.
We now construct tubular neighborhoods
of the deflation sets $\Dd^i_{\Lambda,0}$;
here the codimension is $k = 1$.

\begin{prop}
\label{prop:deflation}
Each $\Dd_{\Lambda,0}^i \subset \ILa$ is a
compact submanifold of codimension $1$ diffeomorphic to $\ILai$,
where $\Lambda_i = \diag(\lambda_1, \ldots, \lambda_{i-1},
\lambda_{i+1}, \ldots \lambda_n)$.
There exists $\epsilon_\tub > 0$ such that
for $\epsilon \in (0,\epsilon_\tub)$:
\begin{enumerate}[(a)]
\item{the connected components
$\Dd_{\Lambda,\epsilon}^i$ of $\Dd_{\Lambda,\epsilon}$ consist of matrices
$T \in \Dd_{\Lambda,\epsilon}$ for which
$|(T)_{n,n} - \lambda_i| < \sqrt{2} \;\epsilon$;}
\item{the map
$\zeta: \Dd^i_{\Lambda,\epsilon} \to
\Dd^i_{\Lambda,0} \times [-\epsilon, \epsilon]$
given by $\zeta(T) = (\Pi_i(T), \bbb(T))$
is a closed tubular neighborhood of $\Dd_{\Lambda,0}$;}
\item{there is a constant $C_\bbb>0$ such that for all
$T \in \Dd^i_{\Lambda,\epsilon}$,
\[ | b(T)| \le \dist(T,\Dd^i_{\Lambda,0}) \le
\| T - \Pi_i(T)\| \le C_\bbb |b(T)|.\]}
\end{enumerate}
\end{prop}

\smallskip
\proof
We first show that the gradient of the restriction
$\bbb|_{\ILa}$ at a point $T_{\Dd} \in \Dd_{\Lambda,0}$ is not zero.
Consider the characteristic polynomial
along the line $T_{\Dd} + tB$:
this is a smooth even function of $t$
and therefore $B$ is tangent to $\ILa$ at $T_{\Dd}$,
the point on which $t=0$.
On the other hand, the directional derivative of $\bbb$
along  the same line equals 1.
Thus $\Dd_{\Lambda,0} \subset \ILa$ is a
submanifold of codimension $1$.
The diffeomorphism with $\ILai$ takes $T$ to
$\hat T$, the leading $(n-1) \times (n-1)$ principal minor of $T$.

Assume $\epsilon < \gamma/(2 \sqrt{2})$.
Consider matrices $T \in \Dd_{\Lambda,\epsilon}$ and $S =T - \bbb(T) B$,
so that $(T)_{n,n}$ is an eigenvalue of $S$.
By the Wielandt-Hoffman theorem,
there exists an index $i$ for which
$|(T)_{n,n} - \lambda_i| < \sqrt{2} \epsilon$,
defining the sets
$\Dd^i_{\Lambda,\epsilon}$
(at this point we do not yet know
that $\Dd^i_{\Lambda,\epsilon}$ is connected).



For $T_{\Dd} \in \Dd^i_{\Lambda,0}$,
the derivative $D\Pi_i(T_{\Dd})$
equals the identity on the subspace tangent
to $\Dd^i_{\Lambda,0}$ and has a kernel of dimension $1$.
Thus, for sufficiently small $\epsilon_\tub$, item (b) holds.
This also proves that each $\Dd^i_{\Lambda,\epsilon}$ is connected,
completing the proof of item (a).

The first two inequalities in (c) are trivial. Now
\[ \| T - \Pi_i(T) \| =
\| \zeta^{-1}(\Pi_i(T),\bbb(T)) -
\zeta^{-1}(\Pi_i(T),0) \| \le
C_\bbb | \bbb(T) |, \]
where the derivative of $\zeta^{-1}(T_{\Dd}, \delta)$
with respect to the second
coordinate is bounded by $C_\bbb$ on the compact set
$\Dd_{\Lambda,0} \times [-\epsilon_\tub, \epsilon_\tub]$.
\qed


The diffeomeorphism $\zeta$ defines \textit{tubular coordinates}
for $T \in \Dd^i_{\Lambda,\epsilon}$:
the matrix $\Pi_i(T) \in \Dd^i_{\Lambda,0} \approx \Tt_{\Lambda_i}$
and $\bbb(T)$.
Under tubular coordinates, $QR$ steps with shift are given by a simple formula.

\begin{coro}
\label{coro:tubularQR}
Consider $\Lambda$, $i$ and $\epsilon \in (0,\epsilon_\tub)$. Then
\begin{align*}
\zeta \circ F_s \circ \zeta^{-1}:
\Dd^i_{\Lambda,0} \times [-\epsilon, \epsilon] &\to
\Dd^i_{\Lambda,0} \times [-\epsilon, \epsilon] \\
(T,b) &\mapsto \left( F_s(T), \frac{(R_\star)_{n,n}}{(R_\star)_{n-1,n-1}} \,b \right)
\end{align*}
where $\zeta^{-1}(T,b) - sI = Q_\star R_\star$.
\end{coro}

\proof
This follows directly from Lemma \ref{lemma:basic} and Propositions \ref{prop:pia} and \ref{prop:deflation}.
\qed

\section{Convergence to deflation}

Sufficiently thin deflation neighborhoods $\Dd_{\Lambda,\epsilon}^i$ are
invariant under $F_s$ for $s \approx \lambda_i$.

\begin{prop}
\label{prop:epsilonF}
Given $C > 0$, there exists
$\epsilon_\inv \in (0, \epsilon_\tub)$,
such that for any $\epsilon \in (0,\epsilon_\inv)$ and
$s \in [\lambda_i - C\,\epsilon,\lambda_i + C\,\epsilon]$ we have
$F_s(\Dd^i_{\Lambda,\epsilon}) \subset \interior(\Dd^i_{\Lambda,\epsilon/2})$.

For a simple shift strategy
$\sigma: \ILa \to \RR$,
there exists $\epsilon_\inv > 0$ such that
if $\epsilon \in (0,\epsilon_\inv)$ then
$F_\sigma(\Dd^i_{\Lambda,\epsilon}) \subset
\interior(\Dd^i_{\Lambda,\epsilon/2})$.
\end{prop}

In particular, $F_s$ is well defined in $\Dd^i_{\Lambda,\epsilon}$ for
$\epsilon \in (0,\epsilon_\inv)$.

\smallskip

\proof
Recall that $F_s(\Dd_{\Lambda,0}^i) = \Dd_{\Lambda,0}^i$.  
From Lemma \ref{lemma:sislambda}, the derivative of $\bbb \circ \Phi_\star$ is zero
at $\Dd_{\Lambda,0}^i \times \{ \lambda_i \}$.
Compactness of $\Dd_{\Lambda,0}^i$ thus implies that
in a sufficiently small neighborhood of
$\Dd_{\Lambda,0}^i \times \{ \lambda_i \}$
we have $|\bbb(F_s(T))| \le |\bbb(T)|/3$.

Now consider a simple shift strategy $\sigma$:
by condition (II) there exists $C_\sigma > 0$ such that
$|\sigma(T) - \lambda_i| < C_\sigma \bbb(T)$;
apply the first statement with $C = C_\sigma$.
\qed


Thus, $F_\sigma$ squeezes neighborhoods
$\Dd^i_{\Lambda,\epsilon}$ at least linearly.
Condition (I) and smoothness imply a stronger version of (II).
We do not want to assume, however,
that $\Dd_{\Lambda,0} \cap \calS_{\sigma} = \emptyset$:
after all, this is not true even for Wilkinson's shift.
We need a more careful statement.

\begin{lemma}
\label{lemma:wlip}
Consider a shift strategy $\sigma$ and
$\epsilon_\inv$ as in Proposition \ref{prop:epsilonF}.
For a compact set $\Kk \subset
\Dd^i_{\Lambda,\epsilon_\inv} \smallsetminus
(\Dd^i_{\Lambda,0} \cap \calS_\sigma)$,
there exists $C_{\Kk}$ such that for all $T \in \Kk$ we have
$|\sigma(T) - \lambda_i| \le C_{\Kk} \bbb(T)^2 $.
\end{lemma}

\proof
Let $\Kk_{\Dd} = \Kk \cap \Dd^i_{\Lambda,0}$;
enlarge $\Kk_{\Dd}$ along $\Dd_{\Lambda,0}^i$
to obtain another compact set
$\Kk_{1} \subset \Dd_{\Lambda,0}^i \smallsetminus \calS_\sigma$,
$\Kk_{\Dd} \subset \interior_{\Dd_{\Lambda,0}^i}(\Kk_{1})$.
Fatten $\Kk_1$ along fibers to define
$\tilde\Kk_1 = \zeta^{-1}(\Kk_1 \times [-\epsilon,\epsilon])$,
$\epsilon \in (0,\epsilon_\inv)$, which, without loss,
still avoids $\calS_\sigma$.
For each $T_\Dd \in \Kk_1$,
consider the function $h_{T_{\Dd}}(b) = \sigma(\zeta^{-1}(T_\Dd,b))$,
obtained by restricting $\sigma$ to a fiber of $\Dd_{\Lambda,\epsilon}^i$.
Each $h_{T_{\Dd}}$ is smooth and even (from condition (I))
and therefore satisfies
$|h_{T_{\Dd}}(b) - \lambda_i| \le C_{T_{\Dd}} |b|^2$.
By compactness, there exists $C_{\Kk_1}$ such that
$|h_{T_{\Dd}}(b) - \lambda_i| \le C_{\Kk_1} |b|^2$ for all $T_{\Dd} \in \Kk_1$.
In other words, there exists $C_{\tilde\Kk_1}$ such that
$|\sigma(T) - \lambda_i| \le C_{\tilde\Kk_1} |\bbb(T)|^2$
for all $T \in \tilde\Kk_1$.
The estimate for $T \notin \tilde\Kk_1$ is trivial.
\qed


\smallskip

{\noindent \bf Proof of Theorem \ref{theo:squeeze}:}
Take $\epsilon = \epsilon_\inv$ as in Proposition \ref{prop:epsilonF}
so that $\Dd^i_{\Lambda,\epsilon}$ is invariant under $F_\sigma$.

Let $\varphi = \bbb \circ \Phi_\star$.
We compute the Taylor expansion of $\varphi(T,s)$ at $(T_{\Dd},\lambda_i)$,
$T_{\Dd} \in \Dd_{\Lambda,0}^i$:
from Lemma \ref{lemma:sislambda}, the gradient of $\varphi$
at $(T_{\Dd},\lambda_i)$ is zero.
Thus, up to a third order remainder,
\begin{align*}
\varphi(T, s) &=
\varphi(T_{\Dd}, \lambda_i) +
\frac{1}{2}\varphi_{T,T}(T_{\Dd}, \lambda_i)(T - T_{\Dd},T - T_{\Dd}) + \\
&\quad + \varphi_{T,s}(T_{\Dd}, \lambda_i)(T - T_{\Dd},s - \lambda_i) +
\frac{1}{2} \varphi_{s,s}(T_{\Dd}, \lambda_i)(s - \lambda_i,s - \lambda_i)
+ \\ &\quad + \hbox{Rem}_3 ( T - T_{\Dd},s - \lambda_i).
\end{align*}
Now, $\varphi(T_{\Dd}, \lambda_i)=0$
and, again from Lemma \ref{lemma:sislambda},
$\varphi(T,\lambda_i) = 0$ for all $T \in \ILa$,
hence $\varphi_{T,T}(T_{\Dd}, \lambda_i)=0$.
Let $C_\sigma$ be the constant in condition (II)
of the definition of a simple shift strategy.
By compactness, there exists $C_1 >0$ such that
for all $T_{\Dd} \in \Dd^i_{\Lambda,0}$,
$T \in \Dd^i_{\Lambda,\epsilon}$ and
$s \in [\lambda_i - C_\sigma\,\epsilon, \lambda_i + C_\sigma\,\epsilon]$,
we have
\[ | \varphi(T,s) | \le
C_1 |s - \lambda_i| ( \| T - T_{\Dd} \| + |s - \lambda_i|) \]
We now apply this estimate for $T_{\Dd} = \Pi_i(T)$,
where $T \in \Dd^i_{\Lambda,\epsilon}$.
By Proposition \ref{prop:deflation}, since $\epsilon < \epsilon_\tub$,
$ \| T - T_{\Dd} \| = \| T - \Pi_i(T)\| \le C_\bbb |\bbb(T)|$
and therefore
\[ | \varphi(T,s) | \le
C_1 |s - \lambda_i| ( C_\bbb |\bbb(T)| + |s - \lambda_i|) \]
implying the quadratic estimate
\[ |\bbb(F_\sigma(T))| =
| \varphi(T,\sigma(T)) | \le
C_1 |\sigma(T) - \lambda_i| ( C_\bbb |\bbb(T)| + |\sigma(T) - \lambda_i|) \le
C_q |\bbb(T)|^2. \]
Using Lemma \ref{lemma:wlip}
instead of condition (II) yields the cubic estimate in (c).
\qed

As a corollary, we obtain the well-known fact that, near deflation,
the rate of convergence of Rayleigh's shift is cubic.
Similarly, the mixed Wilkinson-Rayleigh strategy has cubic convergence.
The rate of convergence for Wilkinson's strategy is far subtler.

\section{Deflationary strategies}

On the way to prove Theorem \ref{theo:connect},
we construct a larger invariant set for $F_\sigma$.
Let $\Uu_\Lambda \subset \ILa$ be the set of unreduced matrices;
for $\epsilon > 0$,
let $\Uu_{\Lambda,\epsilon} =
\Uu_\Lambda \cup \interior(\Dd_{\Lambda,\epsilon})$.
Notice that $\Uu_{\Lambda,\epsilon}$ is open, dense and path-connected.

\begin{lemma}
\label{lemma:ULaeps}
For a shift strategy $\sigma: \ILa \to \RR$,
$\epsilon_\inv$ as in Proposition \ref{prop:epsilonF}
and $\epsilon \in (0,\epsilon_\inv)$,
the open set $\Uu_{\Lambda,\epsilon}$ is invariant under $F_\sigma$.
\end{lemma}

\proof
If $T \in \Uu_\Lambda$ and $\sigma(T)$ is not in the spectrum then
$F_\sigma(T)$ is (well defined and) unreduced.
If $T \in \Uu_\Lambda$ and $\sigma(T) = \lambda_i$ then
$F_\sigma(T) \in \Dd^i_{\Lambda,0} \subset \Uu_{\Lambda,\epsilon}$.
Finally, if $T \in \interior(\Dd^i_{\Lambda,\epsilon})$ then,
by  Proposition \ref{prop:epsilonF},
$F_\sigma(T) \in \interior(\Dd^i_{\Lambda,\epsilon/2}) \subset
\Uu_{\Lambda,\epsilon}$.
\qed

Notice that we do not assume $\sigma$ or $F_\sigma$ to be continuous.

A simple shift strategy $\sigma$ is \textit{deflationary}
if for any $T \in \Uu_{\Lambda,\epsilon_\inv}$
there exists $K \in \NN$  such that
$F_\sigma^K(T) \in \Dd_{\Lambda,\epsilon_\inv}$.

Rayleigh's strategy is known not to be deflationary.
The following well known estimate (\cite{HP} and \cite{Parlett}, section 8-10)
implies that Wilkinson's strategy is not only deflationary but uniformly so,
in the sense that there exists $K$ with
$F_\omega^K(\Uu_{\Lambda,\epsilon_\inv}) \subset \Dd_{\Lambda,\epsilon_\inv}$.  As a corollary, the mixed
Wilkinson-Rayleigh strategy is also uniformly deflationary provided
$\epsilon > 0$ is sufficiently small.


\begin{fact}
\label{fact:parlett}
For $T \in  \Tt$ and $k \in \NN$,
\[ |\bbb(F_\omega^k(T))|^3 \leq
\frac{|\bbb(T)^2 (T)_{n-1,n-2}|}{(\sqrt{2})^{k-1}}.\]
\end{fact}

In \cite{Parlett}, the result is shown for unreduced matrices;
the case $T \in \Uu_{\Lambda,\epsilon_\inv}$
follows by elementary limiting arguments.
Notice that for $T \in \ILa$,
the numerator $|\bbb(T)^2 (T)_{n-1,n-2}|$ is uniformly bounded.


We now prove Theorem \ref{theo:connect}:
if the shift strategy $\sigma: \ILa \to \RR$ is continuous then it is not
deflationary.

{\noindent\bf Proof of Theorem \ref{theo:connect}:}
Fix $\epsilon = \epsilon_\inv/2$.
Let $\Bb^i \subset \Uu_{\Lambda,\epsilon}$ be
the basins of attraction of each invariant neighborhood
$\Dd^i_{\Lambda,\epsilon}$, i.e.,
$T \in \Bb^i$ if there exists $k \in \NN$
such that $F_\sigma^k(T) \in \Dd^i_{\Lambda,\epsilon}$.
The sets $\Bb^i$ are clearly disjoint with
$\Dd^i_{\Lambda,\epsilon} \subset \Bb^i$.
If $\sigma$ is continuous,
they are also open subsets of $\Uu_{\Lambda,\epsilon}$
since $\Bb_i = \bigcup_k F_\sigma^{-k}(\interior(\Dd^i_{\Lambda,\epsilon}))$.
If $\sigma$ is deflationary,
$\bigcup_i \Bb_i = \Uu_{\Lambda,\epsilon}$.
Thus, if $\sigma$ is both continuous and deflationary
then $\Uu_{\Lambda,\epsilon}$ is not connected, a contradiction.
\qed

\section{Dynamics of shifts for a.p.\ free matrices}

From the previous section, cubic convergence may be lost when the orbit
$F_\sigma^k(T)$ passes near the set $\calS_\sigma \cap \Dd_{\Lambda,0}$.
Our next task is to measure when this happens,
by studying the dynamics associated to a shift strategy
in a deflation neighborhood, i.e.,
the iterates of $F_\sigma:\Dd^i_{\Lambda,\epsilon} \to \Dd^i_{\Lambda,\epsilon}$,
$\epsilon \in (0,\epsilon_\inv)$.
Most of what we need can be read in the projection onto $\Dd^i_{\Lambda,0}$,
where $F_\sigma$ coincides with $F_{\lambda_i}$.

A matrix $T \in \Tt$ with simple spectrum is \textit{a.p.\ free}
if no three eigenvalues are in arithmetic progression and a.p.\ otherwise.
Different kinds of spectra lead to different dynamics:
in this section we handle the a.p.\ free case, clearly a generic restriction.
Let $\tilde{T}$ be the leading principal $(n-1) \times (n-1)$ minor of $T$.
The following result is standard.

\begin{prop}
\label{prop:APfreedynamics}
Let $\Lambda \in \Tt$ be an $n \times n$ diagonal a.p.\ free matrix
with spectrum $\lambda_1 < \cdots < \lambda_n$.
For each $i$,
consider $F_{\lambda_i}: \Dd^i_{\Lambda,0} \to \Dd^i_{\Lambda,0}$
as above. For any ${T} \in \Dd^i_{\Lambda,0}$, the sequence
$(F_{\lambda_i}^k(T))$ converges to a diagonal matrix.
\end{prop}

\proof
The map $F_{\lambda_i}$ on $\Dd^i_{\Lambda,0}$ amounts to a $QR$ step
with shift $\lambda_i$ on $\tilde{T}$,
which has eigenvalues $\lambda_j$, $j \ne i$.
The a.p.\ free hypothesis implies that the absolute values of the
eigenvalues of $\tilde T - \lambda_i I$ are distinct.
If $\tilde T$ is unreduced then, as is well known,
the standard $QR$ iteration converges to a diagonal matrix,
with diagonal entries in decreasing order of absolute value.
More generally, if $\tilde T$ is reduced, apply the above result to each
unreduced sub-block.
\qed

We shall use \textit{height functions} for the $QR$ steps $F_s$, $s$ near $\lambda_i$, i.e.,
functions $H_i: \Dd^i_{\Lambda,\epsilon} \to \RR$
with $H_i(F_s(T)) > H_i(T)$ provided $T$ is not diagonal.
Such height functions and related scenarios
have been considered in \cite{BBR}, \cite{DRTW}, \cite{LT} and \cite{Tomei}.

The matrix $W = \diag(w_1,\ldots,w_n)$ is a \textit{weight matrix}
if $w_1 > \cdots > w_n$.
Since $\Lambda$ is a.p.\ free,
there exists $\epsilon_{\textrm{\rm ap}} \in (0,\epsilon_\inv)$ such that
if $s \in \Ii_i = [\lambda_i - \epsilon_{\textrm{\rm ap}},
\lambda_i + \epsilon_{\textrm{\rm ap}}]$ then
the numbers $|\lambda_j - s|$ are distinct and their order does not depend on $s$.

\begin{prop}
\label{prop:height}
Let $\Lambda$ be an a.p.\ free diagonal matrix,
$W$ a weight matrix and $\epsilon_{\textrm{\rm ap}}$ as above.
For $\delta_H > 0$, set $\eta_i(x) = \log((x-\lambda_i)^2 + \delta_H)$
and let $H_i: \Dd^i_{\Lambda,\epsilon_{\textrm{\rm ap}}} \to \RR$
be defined by $H_i(T) = \trace(W \eta_i(T))$.
There exists $\delta_H > 0$ such that
\[ \max_{T \in \partial \Dd^i_{\Lambda,\epsilon_{\textrm{\rm ap}}}} H_i(T)
< \min_{T \in \Dd^i_{\Lambda,0}} H_i(T) \]
and, for any $s \in \Ii_i$,
$H_i$ is a height function
for $F_s: \Dd^i_{\Lambda,\epsilon_{\textrm{\rm ap}}}
\to \Dd^i_{\Lambda,\epsilon_{\textrm{\rm ap}}}$.
\end{prop}

Here,
$\eta_i(T) = X \diag(\eta_i(\lambda_1), \ldots, \eta_i(\lambda_n)) X^{-1}$
for $T = X \Lambda X^{-1}$ so that if $p$ is a polynomial and
$\eta_i(\lambda_j) = p(\lambda_j)$ for $j = 1, \ldots, n$
then $\eta_i(T) = p(T)$.
The only conditions on $\eta_i$ which will be used in the proof are
that $|\lambda_j - \lambda_i| < |\lambda_k - \lambda_i|$ implies
$\eta_i(\lambda_j) < \eta_i(\lambda_k)$
and that $\eta_i(\lambda_i)$ is very negative (for small $\delta_H$).

The proof requires some basic facts about $f$-$Q_\star R_\star$ \textit{steps};
these facts will not be used elsewhere.
For a real diagonal matrix $\Lambda$ with simple spectrum,
let $\Oo_\Lambda$ be the set of
all real symmetric matrices similar to $\Lambda$;
it is well known that $\Oo_\Lambda$ is a smooth compact manifold.
The $f$-$Q_\star R_\star$ \textit{step} applied to a matrix
$S \in \OLa$ is the map $F_f: \Aa_{\Lambda,f} \to \OLa$
defined by $F_f(S) = Q_\star^\ast S Q_\star$,
where $Q_\star$ is obtained from the factorization
$f(S) = Q_\star R_\star$ and $S \in \Aa_{\Lambda,f}$
if and only if $f(S)$ is almost invertible.
If $T \in \ILa \cap \Aa_{\Lambda,f}$ then 
$F_f(T) \in \ILa$ (use the same proof as in Lemma \ref{lemma:basic}).
The maps $F_s: \ILa \to \ILa$ defined above
correspond to restrictions of $F_f$ for $f(x) = x - s$.

%

For a continuous function $h: \RR \to \RR$,
if $S \in \Oo_\Lambda$ then
the matrix function $h(S)$ belongs to $\Oo_\Mu$,
where $\Mu = h(\Lambda)$.
With the obvious abuse of notation,
we have a diffeomorphism $h: \Oo_\Lambda \to \Oo_\Mu$
provided $h$ is injective in the spectrum of $\Lambda$.

\begin{lemma}
\label{lemma:ftilf}
For $h$ injective in the spectrum of $\Lambda$,
consider the diffeomorphism $h: \Oo_\Lambda \to \Oo_\Mu$,
where $\Mu = h(\Lambda)$.
Let $f$ and $\tilde f$ be continuous functions
defined in neighborhoods of the spectra
of $\Lambda$ and $\Mu$, respectively,
satisfying $\tilde f(h(\lambda_j)) = f(\lambda_j)$ for each $j$
with $QR$ steps $F_f:  \Oo_\Lambda \to \Oo_\Lambda$
and $F_{\tilde f}: \Oo_\Mu \to \Oo_\Mu$.
Then $h \circ F_f =  F_{\tilde f} \circ h$.
\end{lemma}

\proof
The hypothesis implies that, for $T \in \Oo_\Lambda$,
$f(T) = \tilde f(h(T)) = QR$ and hence
$F_f(T) = Q^\ast T Q$ and $F_{\tilde f}(h(T)) = Q^\ast h(T) Q$.
Thus $h(F_f(T)) = F_{\tilde f}(h(T))$.
\qed

Let $I_r$ be the $n \times n$ truncated identity matrix,
i.e., $(I_r)_{i,i} = 1$ for $i \le r$, other entries being equal to zero.

\begin{lemma}
\label{lemma:increasing}
Let $\Mu$ be a diagonal matrix with simple spectrum
and $\tilde f: \RR \to \RR$ be a function
for which $\mu_i < \mu_j$ implies
$|\tilde f(\mu_i)| < |\tilde f(\mu_j)|$.
Consider the $\tilde f$-$QR$ step
$F_{\tilde f}: \Aa_{\Mu,\tilde f} \to \Oo_\Mu$.
For any $S \in \Aa_{\Mu,\tilde f}$ and $r = 1, \ldots, n - 1$,
$\trace(I_r F_{\tilde f}(S)) \ge \trace(I_r S)$.
For $r = 1$, equality only holds if $(S)_{1,j} = 0$ for all $j > 1$.
\end{lemma}

This argument follows closely the first proof in \cite{DRTW}.

\smallskip

\proof
Let $V_r$ be the range of $I_r$ and
$\mu_{r,j}(S)$ be the eigenvalues
of the leading principal $r \times r$ minor of $S$,
listed in nondecreasing order.
We claim that $\mu_{r,j}(F_{\tilde f}(S)) \ge \mu_{r,j}(S)$,
which immediately implies $\trace(I_r  F_{\tilde f}(S)) \ge \trace(I_r S)$.
Recall that $F_{\tilde f}(S) = Q_\star^\ast S Q_\star$
where $Q_\star R_\star = \tilde f(S)$.
Let $U$ be an upper triangular matrix such that
$Q_\star u = \tilde f(S) U u$ for $u \in V_r$.
By min-max,
\begin{gather*} \mu_{r,j}(S) =
\max_{ \begin{matrix} \scriptstyle A \subset V_r \\ \scriptstyle \dim(A) = r+1-j \end{matrix}} \;
\min_{u \in A \smallsetminus \{0\}}
\frac{\langle u, S u \rangle}{\langle u, u \rangle}, \\
\mu_{r,j}(F_{\tilde f}(S)) =
\max_{A} \min_{u}
\frac{\langle u, F_{\tilde f}(S) u \rangle}{\langle u, u \rangle} =
\max_{A} \min_{u}
\frac{\langle \tilde f(S)U u, S \tilde f(S)U u \rangle}
{\langle \tilde f(S)U u, \tilde f(S)U u \rangle} \\
=  \max_{A' = U A} \;\min_{u' \in A' \smallsetminus \{0\} }
\frac{\langle \tilde f(S)u', S \tilde f(S)u' \rangle}
{\langle \tilde f(S)u', \tilde f(S)u' \rangle}
\end{gather*}
Notice that since $U$ is upper triangular,
the map taking $A \subset V_r$ to $A' = U A$
is a bijection among subspaces of $V_r$ of given dimension.
Since $S$ and $\tilde f(S)$ are symmetric and commute,
\[ \mu_{r,j}(F_{\tilde f}(S)) = \max_{A} \min_{u}
\frac{\langle u, S g(S) u \rangle}{\langle u, g(S) u \rangle}, \]
where $g(x) = (\tilde f(x))^2$.
The claim now follows from the inequality
\[ \langle u, u \rangle \langle u, S g(S) u \rangle -
\langle u, S u \rangle \langle u, g(S) u \rangle \ge 0. \]
Diagonalize $S= Q^\ast \Mu Q$ and
$g(S) = Q^\ast g(\Lambda)Q$ and write $Qu = (x_1,\ldots,x_n)$ so that
\[ 2 \left(\langle u, u \rangle  \langle u, S g(S) u \rangle -
\langle u, S u \rangle \langle u, g(S) u \rangle \right) =
\sum_{k,\ell} (\mu_k - \mu_\ell)
(g(\mu_k) - g(\mu_\ell)) x_k^2 x_\ell^2 \ge 0.\]
Consider now equality for the case $r = 1$.
Notice that, by hypothesis, if $k \ne \ell$ then
$(\mu_k - \mu_\ell)(g(\mu_k) - g(\mu_\ell)) > 0$.
In the max-min formula for $\trace(I_1 S) = \mu_{1,1}(S)$,
it suffices to take $u = e_1$.
Equality therefore holds only if $Qe_1$ is a canonical vector,
which implies $(S)_{1,j} = 0$ for all $j > 1$.
\qed

{\nobf Proof of Proposition \ref{prop:height}:}
For all $s \in \Ii_i$ and any
distinct eigenvalues $\lambda_j$ and $\lambda_k$,
$|\lambda_j - \lambda_i| < |\lambda_k - \lambda_i|$
if and only if $\eta_i(\lambda_j) < \eta_i(\lambda_k)$.
For $s \in \Ii_i$, $f(x) = x-s$, $h(x) = \eta_i(x)$ and
$\mu_j = \eta_i(\lambda_j)$,
define $\tilde f: \RR \to \RR$ as in Lemma \ref{lemma:ftilf}.
The function $\tilde f$ satisfies the hypothesis
of Lemma \ref{lemma:increasing}:
$\mu_j < \mu_k$ implies $|\tilde f(\mu_j)| < |\tilde f(\mu_k)|$.
Thus, by Lemma \ref{lemma:increasing},
$\trace(W  F_{\tilde f}(S)) \ge \trace(WS)$ for all $S \in \OSa$.
For $T \in \Dd^i_{\Lambda,\epsilon_{\textrm{\rm ap}}}$, take $S = h(T)$:
by Lemma \ref{lemma:ftilf}, $F_{\tilde f}(h(T)) = h(F_f(T))$
and therefore $\trace(W h(F_f(T))) \ge \trace(W h(T))$.
Again by Lemma \ref{lemma:increasing},
equality happens only if $T$ is diagonal.
Thus, $H_i$ is a height function.
Finally, choosing $\delta_H$ sufficiently small
guarantees that $H_i$ is large in $\Dd^i_{\Lambda,0}$
and small in $\partial\Dd^i_{\Lambda,\epsilon_{\textrm{\rm ap}}}$,
completing the proof.
\qed

Thus, simple shift strategies admit height functions near the deflation set.
Our reason for constructing a height function is to control the time the
sequence $(F_\sigma^k(T))$ stays in a compact set.

Assuming $\Lambda$ to be a.p.\ free,
for a shift strategy $\sigma: \ILa \to \RR$
set $\epsilon_\sigma = \epsilon_{\textrm{\rm ap}}/(1+C_\sigma)$
(where $C_\sigma$ is the constant in the definition of a simple shift strategy).
Notice that $T \in  \Dd^i_{\Lambda,\epsilon_{\sigma}}$
implies $\sigma(T) \in \Ii_i =
[\lambda_i - \epsilon_{\textrm{\rm ap}}, \lambda_i + \epsilon_{\textrm{\rm ap}}]$.

\begin{coro}
\label{coro:byebyecompact}
Let $\Lambda$ be a real diagonal $n \times n$ a.p.\ free matrix,
$\sigma$ a simple shift strategy and
$\Dd^i_{\Lambda,\epsilon_{\sigma}}$ as above.
Let $\Kk \subset \Dd^i_{\Lambda,\epsilon_{\sigma}}$ be a compact set
with no diagonal matrices:
there exists $K \in \NN$ such that
for all $T \in \Dd^i_{\Lambda,\epsilon_{\sigma}}$
there are at most $K$ points of the form $F_\sigma^k(T)$ in $\Kk$.
\end{coro}

The plan is to take $\Kk$ containing
$\calS_\sigma \cap \Dd^i_{\Lambda,\epsilon_{\sigma}}$:
the hypothesis in Theorem \ref{theo:big}
that diagonal matrices do not belong to
the singular support $\calS_\sigma$ is then natural.

\proof
Let $m_-$ be the minimum jump in $\Kk$ and
$m_+$ the size of the image of $H_i$:
\[ m_- = \inf_{T \in \Kk,\; s \in \Ii_i} H_i(F_s(T)) - H_i(T), \quad
m_+ = \sup_{T \in \Dd^i_{\Lambda,\epsilon_{\sigma}}} H_i(T) -
\inf_{T \in \Dd^i_{\Lambda,\epsilon_{\sigma}}} H_i(T). \]
By Proposition \ref{prop:height}
and the compactness of $\Kk \times \Ii_i$, $s > 0$:
take $K$ such that $Km_- > m_+$.
For a given $T$, let 
$X = \{ k \in \NN \;|\; F_\sigma^k(T) \in \Kk \}$:
we have
\[ m_+ \ge \sum_{k \in X} H_i(F_\sigma^{k+1}(T)) - H_i(F_\sigma^k(T))
\ge |X| m_- \]
and therefore $|X| < K$.
\qed

{\nobf Proof of Theorem \ref{theo:big}:}
Let $\Kk_1, \Kk_2 \subset \Dd^i_{\Lambda,\epsilon_\sigma}$
be compact sets with
$\Kk_1 \cup \Kk_2 = \Dd^i_{\Lambda,\epsilon_\sigma}$,
$\calS_\sigma \cap \Dd^i_{\Lambda,0}$ disjoint from $\Kk_1$
and with no diagonal matrices in $\Kk_2$.
By Theorem \ref{theo:squeeze},
there exists $C_{\Kk_1} > 0$ such that
$|\bbb(F_\sigma(T))| \le C_{\Kk_1} |\bbb(T)|^3$ for all $T \in \Kk_1$.
By Corollary \ref{coro:byebyecompact},
there exists $K_2 \in \NN$ such that,
given $T \in \Dd^i_{\Lambda,\epsilon_\sigma}$,
at most $K_2$ points of the form $F_\sigma^k(T)$ belong to $\Kk_2$.
In particular, there are at most $K_2$ values of $k$
for which the estimate
$|\bbb(F_\sigma^{k+1}(T))| \le C_{\Kk_1} |\bbb(F_\sigma^k(T))|^3$
does not hold.
\qed

\section{Convergence properties of a.p.\ spectra}

The aim of this section is to prove Theorem \ref{theo:bigg}.
An a.p.\ matrix $T \in \Tt$ with simple spectrum is
\textit{strong a.p.}\ if three consecutive eigenvalues
are in arithmetic progression and \textit{weak a.p.}\ otherwise.

In the a.p.\ free case discussed in the previous sections,
for an initial condition $T \in \Dd^i_{\Lambda,\epsilon}$,
the sequence $F_\sigma^k(T)$ converges to a diagonal matrix;
this follows from the fact that $\sigma(T) \approx \lambda_i$
for $T \in \Dd^i_{\Lambda,\epsilon}$.
For weak a.p.\ spectra, convergence to a diagonal matrix may not occur.

Assume $\Lambda$ to be weak a.p.
Let $\bbb_2(T) = T_{n-1,n-2}$ be the second-last subdiagonal entry;
for consistency, write $\bbb_1(T) = \bbb(T)$.
For any $i$, there exists a unique index
$c(i)$ such that $\lambda_{c(i)}$ is the eigenvalue closest to $\lambda_i$.
As we shall see, if $T \in \Dd^i_{\Lambda,\epsilon}$ then
\[ \lim_{k \to \infty} \bbb_1(F_\sigma^k(T)) = \lim_{k \to \infty} \bbb_2(F_\sigma^k(T)) = 0, \quad
\lim_{k \to \infty} (F_\sigma^k(T))_{n,n} = \lambda_i; \]
furthermore, if $T$ is unreduced then
\[ \lim_{k \to \infty} (F_\sigma^k(T))_{n-1,n-1} = \lambda_{c(i)}. \]

We begin with a technical lemma concerning the dynamics of steps $F_s$.
Item (b) is a variation of the power method argument used
to study the convergence of lower entries under $QR$ steps.

\begin{lemma}
\label{lemma:tube2}
Let $\Mu = \diag(\mu_1, \ldots, \mu_m)$ be a real diagonal
matrix with simple spectrum and $\Tt_\Mu \subset \Tt$
be the manifold of real $m \times m$ tridiagonal matrices similar to $\Mu$.
Let $I \subset \RR$ be a compact interval.
Assume that there exists $j$, $1 \le j \le m$, such that
\[ \mu_j \notin I, \quad
\max_{s \in I} |\mu_j - s| < \min_{k \ne j, s \in I} |\mu_k - s|. \]
Let $\Dd^j_{\Mu,\epsilon} \subset \Tt_\Mu$
be the $j$-th deflation neighborhood.
\begin{enumerate}[(a)]
\item{There exist $\epsilon > 0$ and $C \in (0,1)$ such that
for all $\epsilon' \in (0,\epsilon)$ and $s \in I$ we have
$F_s(\Dd^j_{\Mu,\epsilon'}) \subset \Dd^j_{\Mu,C\epsilon'}$.}
\item{Consider $T_0 \in \Tt_\Mu$ unreduced,
a sequence $(s_k)$ of elements of $I$ and $\epsilon > 0$.
Define $T_{k+1} = F_{s_k}(T_k)$.
Then there exists $k$ such that $T_k \in \Dd^j_{\Mu,\epsilon}$.}
\end{enumerate}
\end{lemma}

This will be used to study $\bbb_2(T)$ for $T \in \Dd^i_{\Lambda,\epsilon}$,
setting $I = [\lambda_i - \epsilon, \lambda_i + \epsilon]$, $j = c(i)$,
$\Mu = \Lambda_i = \diag(\lambda_1, \ldots, \lambda_{i-1},
\lambda_{i+1}, \ldots, \lambda_n)$, 
with the natural identification
between $\Tt_\Mu$ and $\Dd^i_{\Lambda,0}$.

\smallskip

\proof
Let $\tilde C \in (0,1)$ be such that
\[ \max_{s \in I} |\mu_j - s| < \tilde C \;\min_{k \ne j, s \in I} |\mu_k - s|. \]
Write
\[ r(s,T) = \frac{(R_\star)_{m,m}}{(R_\star)_{m-1,m-1}}, \quad
T - sI = Q_\star R_\star. \]
Recall from Lemma \ref{lemma:basic} and
Corollary \ref{coro:tubularQR} that $\bbb(F_s(T)) = r(s,T)\;\bbb(T)$.
We claim that for all $T \in \Dd^j_{\Mu,0}$ and $s \in I$, $|r(s,T)| \le \tilde C$.
Since $T \in \Dd^j_{\Mu,0}$, $|(R_\star)_{m,m}| = |\mu_j - s|$.
Let $R_-$ be the leading principal minor of $R_\star$ of order $m-1$:
its singular values are $|\mu_k - s|$, $k \ne s$.
In particular, all singular values
are larger that $|(R_\star)_{m,m}|/{\tilde C}$.
Thus
\[ |(R_\star)_{m-1,m-1}| = \| e_{m-1}^\ast R_- \| \ge
\frac{|(R_\star)_{m,m}|}{\tilde C} \| e_{m-1} \| =
\frac{|(R_\star)_{m,m}|}{\tilde C}, \]
proving our claim.
Take  $C = (1+\tilde C)/2$:
by continuity, for sufficiently small $\epsilon > 0$, we have $|r(s,T)| < C$
for all $T \in \Dd^j_{\Mu,\epsilon}$, $s \in I$.
Thus, for $T \in \Dd^j_{\Mu,\epsilon}$ and $s \in I$,
$|\bbb(F_s(T))| \le C\,|\bbb(T)|$;
item (a) follows.


For item (b), write $T_{k+1} = Q_k^\ast T_k Q_k$
where $T_k - s_k I = Q_k R_k$ is a $Q_\star R_\star$ decomposition.
Notice that, by hypothesis, $I$ is disjoint from the spectrum so that
$T_0 - s_0 I$ is invertible.
We have $(T_0 - s_0 I)^{-1} = R^{-1} Q_0^\star$
so the rows of $Q_0^\star$ are obtained from those of $(T_0 - s_0 I)^{-1}$
by Gram-Schmidt from bottom to top.
In particular, $Q_0 e_m = c_0 (T_0 - s_0 I)^{-1} e_m$,
$c_0 > 0$.
More generally, we claim that
\begin{gather*}
P_k e_m =
c (T_0 - s_{k-1} I)^{-1} \cdots (T_0 - s_1 I)^{-1} (T_0 - s_0 I)^{-1} e_m, \\
c > 0, \quad P_k = Q_0 Q_1 \cdots Q_{k-1} \in SO(m).
\end{gather*}
Indeed, by induction and using that $T_1 = Q_0^\ast T_0 Q_0$,
\begin{align*}
P_k e_m &= c' Q_0 (T_1 - s_{k-1} I)^{-1} \cdots (T_1 - s_1 I)^{-1} e_m  \\
&= c' (T_0 - s_{k-1} I)^{-1} \cdots (T_0 - s_1 I)^{-1} Q_0 e_m \\
&= c (T_0 - s_{k-1} I)^{-1} \cdots (T_0 - s_1 I)^{-1} (T_0 - s_0 I)^{-1} e_m.
\end{align*}
For $\alpha = 1, \ldots, m$,
let $v_\alpha$ be the unit eigenvector associated to $\mu_\alpha$.
We claim that
\[ \lim_{k \to \infty} P_k e_m = \pm v_j. \]
Indeed, write $e_m = \sum_{\alpha=1}^m a_\alpha v_\alpha$,
where $a_\alpha = \langle v_\alpha, e_m \rangle$ is the last
coordinate of $v_\alpha$.
It is well known that the last coordinates of the eigenvectors $v_\alpha$
of the unreduced matrix $T$ are nonzero:
in particular, $a_j \ne 0$; assume without loss $a_j > 0$.
We have
\begin{align*}
P_k e_m &=
c (T_0 - s_{k-1} I)^{-1} \cdots (T_0 - s_1 I)^{-1} (T_0 - s_0 I)^{-1} e_m \\
&= c \sum_{\alpha = 1}^m
\frac{a_\alpha}{(\mu_\alpha - s_{k-1})\cdots(\mu_\alpha - s_{0})}
v_\alpha
= c_k \left( v_j + \sum_{\alpha \ne j} b_{k,\alpha} v_\alpha \right),\\
&\qquad c_k > 0,  \quad
b_{k,\alpha} =
\frac{a_\alpha}{a_j}\;\frac{\mu_j - s_{k-1}}{\mu_\alpha - s_{k-1}}\cdots
\frac{\mu_j - s_{0}}{\mu_\alpha - s_{0}}.
\end{align*}
Since $| \mu_j - s_{k-1} |/|\mu_\alpha - s_{k-1}| < \tilde C$
we have $|b_{k,\alpha}| \le (\tilde C)^k\,|a_\alpha/a_j|$ and
therefore $\lim_{k \to \infty} b_{k,\alpha} = 0$,
proving the claim. We have
\begin{gather*}
\lim_{k \to \infty} \bbb(T_k) = \lim_{k \to \infty} (T_k)_{m,m-1} =
\lim_{k \to \infty} e_{n-1}^\ast T_k e_m =
\lim_{k \to \infty} (P_k e_{m-1})^\ast T_0 (P_k e_m)= \\
= \lim_{k \to \infty} (P_k e_{m-1})^\ast \mu_j (P_k e_m) +
\lim_{k \to \infty} (P_k e_{m-1})^\ast (T_0 - \mu_j I) (P_k e_m).
\end{gather*}
The first limit in the last expression is zero because
$P_k e_{m-1} \perp P_k e_{m}$;
the second is zero because $P_k e_{m-1}$ is bounded
and
\[ \lim_{k \to \infty} (T_0 - \mu_j I) (P_k e_m)
= (T_0 - \mu_j I) \lim_{k \to \infty}  (P_k e_m)
= (T_0 - \mu_j I) v_j = 0. \]
\qed

Consider the \textit{double deflation set}
$\Cc_{\Lambda,0} \subset \Dd_{\Lambda,0} \subset \Tt_\Lambda$:
\[ \Cc_{\Lambda,0} = \{ T \in \Tt_{\Lambda} \;|\; \bbb_1(T) =
\bbb_2(T) = 0 \}. \]
For Wilkinson's strategy $\omega$, it turns out that the set $\Cc_{\Lambda,0}$
is disjoint from the singular support $\calS_\omega$.  More generally, if a
shift strategy $\sigma$ satisfies $\Cc_{\Lambda,0} \cap \calS_\sigma =
\emptyset$ then cubic convergence of $F_\sigma$ holds even for weak a.p.\
spectra: this is Theorem \ref{theo:bigg}, which we prove below.

In \cite{LST2}, we show examples of
unreduced tridiagonal $3 \times 3$ matrices
with spectrum $-1, 0, 1$ for which Wilkinson's shift $F_\omega$ converges
quadratically to a reduced but not diagonal matrix
in the singular support $\calS_\omega$.
Similarly, we conjecture that
for strong a.p.\ diagonal $n \times n$ matrices $\Lambda$
there exists a set $\Xx \subset \ILa$ of Hausdorff codimension 1
of unreduced matrices $T$ for which $F_\omega^k(T)$ converges quadratically
to a matrix in $\calS_\omega \cap \Dd_{\Lambda,0}$ with $T_{n-1,n-2} \ne 0$.

%
%

With the natural identification between $\Dd^i_{\Lambda,0}$ and $\Tt_{\Lambda_i}$,
we may consider $\Dd^j_{\Lambda_i,\epsilon_2}$ to be a subset of $\Dd^i_{\Lambda,0}$.
Let
\[ \Cc^{j,i}_{\Lambda,\epsilon_2,\epsilon_1} =
\Dd^i_{\Lambda,\epsilon_1} \cap
\Pi_i^{-1}(\Dd^{j}_{\Lambda_i,\epsilon_2}). \]
For small $\epsilon_1, \epsilon_2 > 0$,
$T \in \Cc^{j,i}_{\Lambda,\epsilon_2,\epsilon_1}$
implies
\[ T_{n-1,n-1} \approx \lambda_j, \quad T_{n,n} \approx \lambda_i,
\quad \bbb_1(T) \le \epsilon_1, \quad \bbb_2(T) \approx 0. \]
These compact sets turn out to be manifolds with corners
but we shall neither prove nor use this fact.
Lemma \ref{lemma:tube2} can be rephrased in terms of the sets $\Cc^{j,i}_{\Lambda,\epsilon_2,\epsilon_1}$.

\begin{coro}
\label{coro:tube3}
Let $\Lambda$ to be weak a.p. spectrum
and $\sigma$ be a simple shift strategy.
There exists $\epsilon > 0$ such that,
for all $i$ and for all $\epsilon_1 \in (0,\epsilon)$:
\begin{enumerate}[(a)]
\item{there exists $C \in (0,1)$ such that,
for all sufficiently small $\epsilon_2 > 0$
we have $F_\sigma(\Cc^{c(i),i}_{\Lambda,\epsilon_2,\epsilon_1}) \subset
\Cc^{c(i),i}_{\Lambda,C \epsilon_2,\epsilon_1}$;}
\item{for all unreduced $T \in \Dd^i_{\Lambda,\epsilon}$
and for all $\epsilon_1, \epsilon_2 > 0$ there exists $k$ such that
$F_\sigma^k(T) \in \Cc^{c(i),i}_{\Lambda,\epsilon_2,\epsilon_1}$.}
\end{enumerate}
\end{coro}

\proof
Combine Lemma \ref{lemma:tube2} with
$\Pi_i \circ F_s = F_s \circ \Pi_i$
(Proposition \ref{prop:pia}).
\qed

{\nobf Proof of Theorem \ref{theo:bigg}:}
From the hypothesis that $\Cc_{\Lambda,0}$ and $\calS_\sigma$ are disjoint
it follows that, for sufficiently small $\epsilon_1, \epsilon_2 > 0$,
the shift strategy $\sigma$ is smooth in $\Cc^{c(i),i}_{\Lambda,\epsilon_2,\epsilon_1}$.
As in Lemma \ref{lemma:wlip}, from a Taylor expansion around $T_0 \in \Dd^i_{\Lambda,0}$,
there exists $C_2$ such that $|\sigma(T)| \le C_2 |\bbb_1(T)|^2$ for all
$T \in \Cc^{c(i),i}_{\Lambda,\epsilon_2,\epsilon_1}$.
As in the proof of Theorem \ref{theo:squeeze}, there exists $C_3$
such that $|\bbb_1(F_\sigma(T))| \le C_3 |\bbb_1(T)|^3$
for all $T \in \Cc^{c(i),i}_{\Lambda,\epsilon_2,\epsilon_1}$.
From item (a) of Corollary \ref{coro:tube3},
$\Cc^{c(i),i}_{\Lambda,\epsilon_2,\epsilon_1}$
is invariant under $F_\sigma$;
from item (b), for all unreduced $T \in \Dd^i_{\Lambda,\epsilon}$
(where $\epsilon$ is sufficiently small)
there exists $K$ such that, for all $k > K$, $F_\sigma^k(T) \in \Cc^{c(i),i}_{\Lambda,\epsilon_2,\epsilon_1}$,
completing the proof.
\qed

\section{Two counterexamples}

In this section we present two examples which show that
natural strengthenings of Theorems \ref{theo:big} and \ref{theo:bigg}
do not hold for Wilkinson's strategy $\omega$.

\begin{figure}[ht]
\begin{center}
\psfrag{Yy}{$\calS_\omega$}
\psfrag{S}{$S$}
\epsfig{height=28mm,file=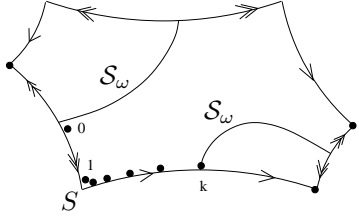}
\end{center}
\caption{\capsize We may have $F_\omega^k(T) \in \calS_\omega$ for large values of $k$.}
\label{fig:notsocubic}
\end{figure}

We use the notation of Section 3.
In Figure \ref{fig:notsocubic}, where $\Lambda = \diag(1,2,4)$,
we indicate a sequence $\tilde F_\omega^k(T)$ 
which enters the deflation neighborhood $\Dd^i_{\Lambda,\epsilon}$
near one diagonal matrix but travels within the neighborhood
towards another diagonal matrix.
Theorem \ref{theo:squeeze} guarantees the cubic decay of the $(3,2)$
entry whenever $\tilde F_\omega^k(T)$ stays away from the singular support $\calS_\omega$.
Consistently with Theorem \ref{theo:big},
this happens for practically all values of $k$.
Notice however that no uniform bound exists on the number of iterations
needed to reach (a neighborhood of) $\calS_\omega$.
As proved in \cite{LST2}, in this instance cubic decay does not hold.
More precisely,
it is {\it not} true that given an a.p.\ free matrix $\Lambda$
there exist $C > 0$ and $K$ such that
$|\bbb(F_\omega^{k+1}(T))| \le C |\bbb(F_\omega^k(T))|^3$ for all $k > K$.

Consider now the weak a.p. spectrum $\Lambda = \diag(-1,0,0.3,1)$ and
\[ T_0 = \begin{pmatrix} 0.3 & 0 \\ 0 & S_0 \end{pmatrix} \in \ILa \]
where $S_0 \in \Tt_{\Lambda_3}$, $\Lambda_3 = \diag(-1,0,1)$,
is an example of unreduced matrix obtained in \cite{LST2}
for which convergence is strictly quadratic, i.e.,
\[ C_- |\bbb(F_\omega^k(S_0))|^2 < |\bbb(F_\omega^{k+1}(S_0))| <
C_+ |\bbb(F_\omega^k(S_0))|^2, \]
for all $k$, where $0 < C_- < C_+$.
Trivially, the analogous estimate holds for $\bbb(F_\omega^k(T_0))$.
By sheer continuity, given $K$, there exists $\epsilon > 0$ such that if
$T \in \ILa$ satisfies $\| T - T_0 \| < \epsilon$  then
\[ C_- |\bbb(F_\omega^k(T))|^2 < |\bbb(F_\omega^{k+1}(T))| <
C_+ |\bbb(F_\omega^k(T))|^2 \]
still holds for all $k < K$.
Thus, the uniform estimate in Theorem \ref{theo:big}
fails for weak a.p. spectra, even for unreduced matrices.

\bigskip\bigskip\bigbreak

{

\parindent=0pt
\parskip=0pt
\obeylines

Ricardo S. Leite, Departamento de Matemática, UFES
Av. Fernando Ferrari, 514, Vitória, ES 29075-910, Brazil

\smallskip

Nicolau C. Saldanha and Carlos Tomei, Departamento de Matem\'atica, PUC-Rio
R. Marqu\^es de S. Vicente 225, Rio de Janeiro, RJ 22453-900, Brazil

\smallskip

rsleite@cce.ufes.br
saldanha@puc-rio.br; http://www.mat.puc-rio.br/$\sim$nicolau/
tomei@mat.puc-rio.br

}

\end{document}